\documentclass[letterpaper,12pt]{siamart190516}

\title{A Quasi-Optimal Spectral Solver for the Heat and Poisson Equations in a Closed Cylinder}
\author{$ $\\David Darrow\thanks{Massachusetts Institute of Technology, Department of Mathematics}\\{\small Advisors: Alex Townsend\thanks{Cornell University, Department of Mathematics}, Grady Wright\thanks{Boise State University, Department of Mathematics}}\\$ $}
\usepackage{mathtools, amsfonts, amssymb, graphicx, hyperref,float,subcaption,fancyhdr,titling,tensor,mathrsfs,xcolor,tikz,tkz-euclide, booktabs}
\usepackage[inline]{enumitem}
\usetikzlibrary{decorations.pathmorphing}

\usepackage[onehalfspacing,nodisplayskipstretch]{setspace}

\newtheorem{algo}{Algorithm}

\newcommand{\vc}[1]{\mathbf{#1}}

\begin{document}

\maketitle

\begin{abstract}
	We develop a spectral method to solve the heat equation in a closed cylinder, achieving a quasi-optimal $\mathcal{O}(N\log N)$ complexity and high-order, \emph{spectral} accuracy. The algorithm relies on a Chebyshev--Chebyshev--Fourier (CCF) discretization of the cylinder, which is easily implemented and decouples the heat equation into a collection of smaller, sparse Sylvester equations. In turn, each of these equations is solved using the alternating direction implicit (ADI) method in quasi-optimal time; overall, this represents an improvement in the heat equation solver from $\mathcal{O}(N^{4/3})$ (in previous Chebyshev-based methods) to $\mathcal{O}(N\log N)$. While Legendre-based methods have recently been developed to achieve similar computation times, our Chebyshev discretization allows for far faster coefficient transforms; we demonstrate the application of this by outlining a spectral method to solve the incompressible Navier--Stokes equations in the cylinder in quasi-optimal time. Lastly, we provide numerical simulations of the heat equation, demonstrating significant speed-ups over traditional spectral collocation methods and finite difference methods.
\end{abstract}
\pagebreak

\section{Introduction}
The forced heat equation models the diffusion of heat through a homogeneous, isotropic medium; it is given by
\begin{equation}\label{eq:heat}
	\big(\partial_t - \alpha\nabla^2\big) T = g(x,t),
\end{equation}
where $T$ is the local temperature field, $\alpha$ is the thermal diffusivity of the medium, and $g$ is the rate of heat generation. More than just heat diffusion, however, the heat equation provides an archetype for diffusion problems in general contexts---for one, viscous or slow-moving (incompressible) fluids obey the \emph{Stokes} equations, a diffusive approximation to the general Navier--Stokes equations when viscous forces dominate the flow:
\begin{equation*}
	\nabla\cdot \vec{v} = 0,\qquad\big(\partial_t - \tfrac{1}{Re}\nabla^2\big) \vec{v} = -\vec{v}\cdot\nabla\vec{v}-\rho^{-1}\nabla p\approx -\rho^{-1}\nabla p.
\end{equation*}
Here, $\vec{v}$ is the fluid velocity, $p/\rho$ is the normalized fluid pressure, and $Re$ is the \emph{Reynolds number} of the flow, a dimensionless group characterizing the relative effects of advective and viscous forces; in making the Stokes approximation, we are assuming $Re\ll 1$. 

As such, the heat equation has seen a great deal of attention in both numerical analysis and in industry applications. For one, finite difference (FD) methods for the heat equation are well-catalogued in textbooks (see \cite{cooper2012introduction} for a modern introduction), and for good reason; FD approaches are \emph{spatially local}, meaning that (a) they can be more easily adapted to new geometries and (b) they allow for equation coefficients (such as the thermal diffusivity) to vary in space. For more recent examples, Berrone has adapted the \emph{finite element} (FE) method to systems with discontinuous coefficients \cite{doi:10.1137/080737058}, and Causley et al. have introduced a powerful convolution-based approach to solve the heat equation and certain nonlinear analogues in optimal time \cite{doi:10.1137/15M1035094}.

Though these approaches allow for fast solvers and great flexibility in the form of the equation and domain, they trade for this by giving only \emph{polynomial} accuracy in the solution itself---that is, the solution generally carries an error of $\mathcal{O}(N^{-M})$, where $N$ is the discretization lattice size and $M$ is a small positive integer \cite{10.5555/1355322}. To achieve better resolution, high-order \emph{spectral collocation methods} are often applied; these approximate the solution space using degree-$N$ polynomials, allowing for far better, $\mathcal{O}(\rho^{-N})$ errors (for some $\rho>1$), but traditionally at the cost of more expensive runtimes \cite{10.5555/3384673}. 

Spectral \emph{coefficient} methods can circumvent these runtimes in highly symmetric domains---namely, when derivatives can be expressed as sparse matrices in a chosen polynomial basis. For one, this means they are particularly suited for periodic domains, where the standard Fourier basis diagonalizes differentiation. If we approximate a function $f$ as $f = \sum f_\vc{k} e^{i\vc{k}\cdot\vc{x}}$, the Laplacian operator acts as $\nabla^2:f_\vc{k}\mapsto -\vc{k}^2f_\vc{k}$, decomposing (\ref{eq:heat}) into a collection of easily-solved scalar equations
\[(\partial_t + \alpha\vc{k}^2)T_\vc{k} = g_\vc{k}.\]
In general, a direction of periodicity allows the $N$-dimensional heat equation to decouple into a collection of $(N-1)$-dimensional equations. 

A critical example of this is given by the 2-dimensional disk, for which a variety of spectral coefficient methods have been developed over the years \cite{PhysRevLett.26.1100,GODON1993171,doi:10.1137/0916061,MATSUSHIMA1995365,verkley_97,torres1999,SHEN2000183,KWAN2009170,doi:10.1137/16M1070207,VASIL201653}. Much like in the periodic case, we can approximate a function as $f=\sum f_{k}(r)e^{ik\theta}$, and the heat equation decomposes as
\begin{equation}\label{eq:disk}
	\left(\partial_t - \alpha(\partial_r^2 + r^{-1}\partial_r +r^{-2}k^2)\right)T_k(r) = g_k(r).
\end{equation}
In turn, the radial function $T_k(r)$ can be expanded in a basis of orthogonal polynomials, such as that of Zernike, Bessel, or Chebyshev; a thorough comparison of bases is given in \cite{BOYD20111408}. Though the radial derivative is not diagonalized, an appropriate choice of basis can reduce (\ref{eq:disk}) to a sparse system of linear equations. More recently, Olver and Townsend introduced a spectral method for the disk (and sphere) that works even with spatially-varying coefficients \cite{doi:10.1137/120865458, doi:10.1137/16M1070207}, and later adapted this method to the unit ball in 3-D \cite{Boull__2020}. 

The general problem of 3-D domains is somewhat more subtle than the disk, as the equations corresponding to (\ref{eq:disk}) now involve at least two independent coordinates. Though spectral (and joint spectral-finite-difference) methods have been developed in the cylinder \cite{gottlieb1977numerical, TAN198581, PULICANI199193, PRIYMAK1995366, doi:10.1137/S1064827595295301,CiCP-5-426} (and see \cite{https://doi.org/10.48550/arxiv.2111.04585} for a recent method developed for the cube), most sacrifice either the spectral accuracy or quasi-optimal computation time possible in the disk. A notable exception is \cite{doi:10.1137/S1064827595295301}, which makes use of the generalized cyclic reduction method of \cite{10.2307/2156231} to achieve both; this requires the use of a Legendre or Chebyshev--Legendre method, however, which has the disadvantage of more expensive value-coefficient transforms than a pure Chebyshev approach \cite{shen_2001}. Their algorithm also requires the introduction of boundary conditions at the central axis of the cylinder, which we avoid in our own discretization. It is worth mentioning that the quasi-optimal method of \cite{VASIL201653} for the disk has been applied to the infinite cylinder in Dedalus \cite{Burns_2020}: a demonstration is shown \href{https://dedalus-project.readthedocs.io/en/v2_master/notebooks/TaylorCouetteFlow.html}{here}. Note, however, that this cylinder is periodic in the $z$-direction, so sidesteps the geometric issues discussed above.

In this paper, we adapt the methods of Olver and Townsend \cite{doi:10.1137/120865458} to the closed cylinder, achieving a quasi-optimal solver for the heat equation with spectral accuracy. Our method is based off of a Chebyshev--Chebyshev--Fourier decomposition ``doubled'' in the $r$-direction, a discretization introduced in \cite{PhysRevLett.26.1100} in the case of the disk that allows $r$ to range from $-1$ to $1$. This discretization carries two significant benefits: it avoids a false boundary along the central axis of the cylinder, and it allows the use of fast Chebyshev transforms to recover function values from polynomial coefficients. It should be noted that this discretization has been observed to cause stability problems at high resolutions (and particularly with nonlinear equations), but these can be mitigated with parity restrictions as discussed in \cite{BOYD20111408} or by imposing partial regularity of the solutions, as discussed in \cite{https://doi.org/10.48550/arxiv.1710.11259}.

As discussed above, for each Fourier mode in our discretization, we must solve an equation of the form (\ref{eq:disk}), though with two independent coordinates, $r$ and $z$: in our framework, this takes the form of a sparse Sylvester equation rather than a system of linear equations. The primary advancement of the present work is the application of the \emph{alternating direction implicit} (ADI) method of \cite{doi:10.1137/0103003} to solve these Sylvester equations. This maintains spectral accuracy but results in an overall $\mathcal{O}(N\log N)$ computation time---this compares favorably with the $\mathcal{O}(N^{7/3})$ necessary to solve these equations directly, or the $\mathcal{O}(N^{4/3})$ Chebyshev-based method given in \cite{10.2307/2156231}.

Finally, we demonstrate a potential application of this method to solve the incompressible Navier--Stokes equations in the cylinder, with the same time and accuracy as with the heat equation. Turbulence modeling highlights the benefits of our approach; the increased digits of precision with spectral accuracy---often 5--10 or more---are essential for turbulence resolution \cite{ZHIYIN201511}.

Our domain of interest for this paper is the closed cylinder in $\mathbb{R}^3$, though the method can be applied equally well to, e.g., the square or the unit ball. The cylinder has many applications throughout the sciences, including studying blood flow, storing fuel in rockets, and understanding properties of hurricanes. Accordingly, the cylindrical coordinates $(r,z,\theta)$ are used throughout. Note that this geometry differs from studying a flow across a cylinder, a problem common in aerodynamics.

Our code was originally developed as part of MIT PRIMES 2017, and is freely available on GitHub \cite{mycode}; since then, these methods have been adapted for fast simulations of active fluids in a ball \cite{https://doi.org/10.48550/arxiv.2103.16638} and to develop improved Poisson solvers in various geometries \cite{10.1093/imanum/drz034}.

\section*{Acknowledgments}
I would like to thank Professors Alex Townsend (Cornell) and Grady Wright (Boise State University), who suggested the original problem and mentored me throughout the research process. I would also like to thank the MIT PRIMES research program, and particularly Professor Pavel Etingof, Dr.~Tanya Khovanova, and Dr.~Slava Gerovitch for making this research possible.

\section{Spectral Approximation of Functions in the Cylinder}\label{sec:approx}
In order to approximate solutions to the heat equation in the cylinder, we first discretize the system by restricting the geometry to $N$ selected \emph{discretization points} within the cylinder. In turn, this allows us to limit the heat equation to a $N$-dimensional space of orthogonal polynomials on the cylinder. This approach, a \emph{spectral method}, provides one unique advantage over more traditional methods: so-called \emph{spectral accuracy}.

While standard simulation methods (such as finite difference and finite element methods mentioned in the preceding section) achieve algebraic accuracy---errors of $\mathcal{O}(N^{-k})$ for some $k>0$ \cite{10.5555/1355322}---spectral accuracy achieves errors of $\mathcal{O}(\rho^{-N})$ for some $\rho>1$ whenever the solution is smooth. Even when $N$ is as small as $1000$, this often represents an improvement over algebraically accurate methods of several orders of magnitude. This accuracy is also uniform throughout the cylinder; some current methods, such as finite element methods, often underresolve the boundary of the region, potentially missing important effects of \textquotedblleft wall-bounded turbulence\textquotedblright{} \cite{Krauthammer1979AccuracyOT}.

Let $\mathit{C}\subset\mathbb{R}^3$ be the closed cylinder written in cylindrical coordinates $(r,z,\theta)$, where $-1\leq r\leq 1$, $-\pi\leq \theta \leq \pi$, and $-1\leq z\leq1$. Let $m,n,p\in \mathbb{N}$ be the number of discretization points in the radial, vertical, and angular directions, respectively, so the total number of discretization points is $N=mnp$.

We discretize using the ``doubled'' Chebyshev--Chebyshev--Fourier grid, introduced in the case of the disk in \cite{PhysRevLett.26.1100} and adapted for the cylinder here:

\begin{definition}\label{grid}
	The \emph{doubled $m\times n\times p$ Chebyshev--Chebyshev--Fourier grid}, denoted by CCF$_{(m,n,p)}$, is the set (in cylindrical coordinates $(r,z,\theta)$) of points $\mathbf{x}_{(j,k,\ell)}$, where
	$$\mathbf{x}_{(j,k,\ell)}=\left(\cos\left(\frac{(m-j-1)\pi}{m-1}\right), \cos\left(\frac{(n-k-1)\pi}{n-1}\right),\frac{(2\ell - p)\pi}{p} \right),$$
	$$j,k,\ell\in\mathbb{Z}, \quad 0\leq j<m, \quad 0\leq k<n, \quad 0\leq \ell<p.$$
\end{definition}

\begin{figure}
	\begin{subfigure}[b]{.33\textwidth}
		\centering
		\includegraphics[scale=.16]{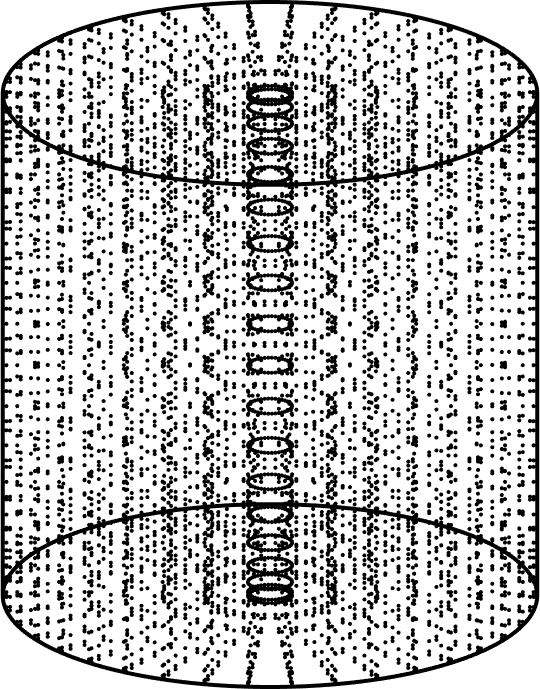}
	\end{subfigure}%
	\begin{subfigure}[b]{.33\textwidth}
		\centering
		\includegraphics[scale=.25]{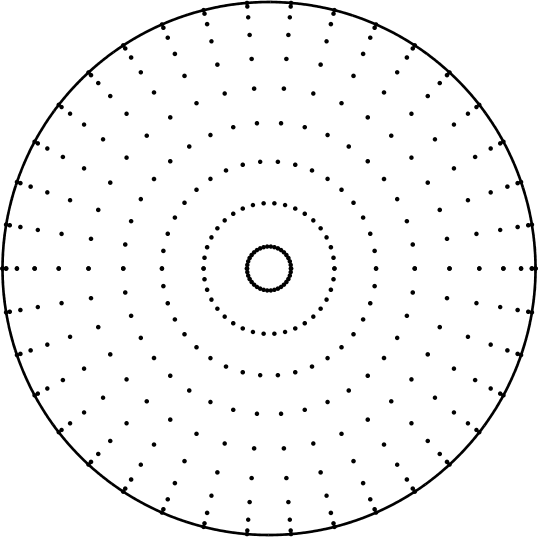}
	\end{subfigure}%
	\begin{subfigure}[b]{.33\textwidth}
		\centering
		\includegraphics[scale=.25]{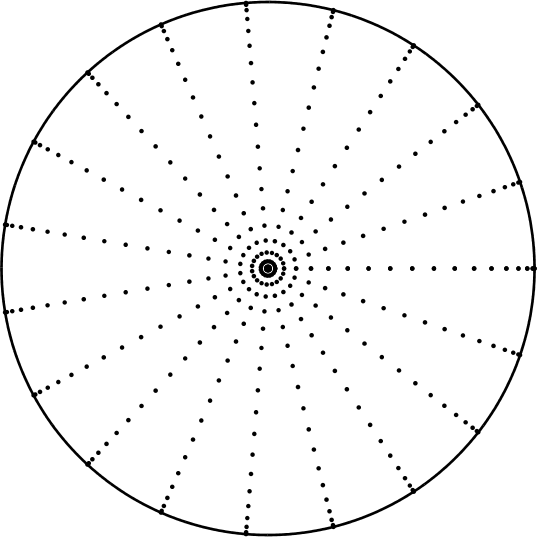}
	\end{subfigure}
\caption{Illustration of the CCF$_{(m,n,p)}$ discretization grid (along with a cross-section), for $m=n=p=20$, contrasted with the corresponding undoubled grid (the current paradigm). Unlike the ``undoubled'' grid on the right, CCF$_{(m,n,p)}$ does not cluster the majority of the points near the centerline of the cylinder.}\label{gridfig}
\end{figure}

This ``doubled'' grid, shown in Figure \ref{gridfig}, can be contrasted with the non-doubled Chebyshev--Chebyshev--Fourier grid in the cylinder, which chooses the radial Chebyshev points between $0$ and $1$. This clusters discretization points near the azimuth (centerline) of the cylinder: a full third of the points are within the center area $1/16^{th}$ the cross sectional size of the region, slowing down computation significantly without providing an increase in approximation accuracy. The doubled grid, however, chooses the radial Chebyshev points between $-1$ and $1$. This grid is an improvement over the standard one; it distributes points almost evenly with respect to cross sectional area. For further information on the CCF grid, see \cite{doi:10.1137/16M1070207}, where it is introduced for spherical and polar geometries. 

We should note that this discretization has been observed to cause stability problems at high resolutions, as it does  (and particularly with nonlinear equations), but these can be mitigated with parity restrictions as discussed in \cite{BOYD20111408}.

Now, we must choose some $N$-dimensional basis to span the set of functions restricted to our discretization grid. As it turns out, the set of ``CCF$_{(m,n,p)}$ polynomials'', defined below, does the trick; the set satisfies necessary existence and uniqueness theorems, allows easy representation of certain differential operators, and interpolates with high accuracy \cite{10.5555/3384673}.

\begin{definition}
	A \emph{CCF$_{(m,n,p)}$ polynomial} is a function $q_{(j,k,\ell)}:\mathit{C}\rightarrow\mathbb{R}^3$ of the following form:
	$$ 
		q_{(j,k,\ell)}(r,z,\theta)=T_{j}(r)T_{k}(z)e^{i (\ell - \lfloor p/2\rfloor) \theta},
	$$
	where $T_K$ is the degree $K$ Chebyshev polynomial of the first kind and $j,k,\ell$ have the same bounds as in Definition \ref{grid}. The complex vector space spanned by these CCF$_{(m,n,p)}$ polynomials is denoted $P_{(m,n,p)}$. 
\end{definition}

Also of use are these two sorts of polynomials: $U_K=C^{(1)}_K$ indicates the degree $K$ Chebyshev polynomial of the second kind, and $C^{(2)}_K$ indicates the degree $K$ ultraspherical polynomial such that the set $\{C^{(2)}_i\}_{i=0}^\infty$ is orthogonal with respect to the weight $(1-x^2)^{3/2}$ on $[-1,1]$. We use the normalization of $C^{(j)}_K$ given in \cite{abramowitz1965handbook}; importantly, we have
\[\partial_xT_K(x) = KU_{K-1}(x),\qquad \partial_x^2T_K(x^2) = 2KC^{(2)}_{K-2}(x).\]
For more information (and explicit formulas) on ultraspherical polynomials, please consult \cite{abramowitz1965handbook}.

Any continuous function $f:\mathit{C}\rightarrow\mathbb{R}$ has a unique interpolant $f^*\in P_{(m,n,p)}$---i.e. $f=f^*$ on CCF$_{(m,n,p)}$---given by
$$
f^*(r,z,\theta) = \sum\limits_{j=0}^{m-1}\sum\limits_{k=0}^{n-1}\sum\limits_{\ell=0}^{p-1} f^{(j,k,\ell)}T_j(r)T_k(z)e^{i(\ell - \lfloor p/2\rfloor)\theta} = \sum f^{(j,k,\ell)}q_{(j,k,\ell)}(r,z,\theta).
$$

Starting with the set of values of $f$ on the grid, we can find the coefficients $f^{(j,k,\ell)}$ through the following procedure. First set $f_{(j,k,\ell)}=f(\mathbf{x}_{(j,k,\ell)})$, and perform the discrete cosine transform of type I (DCT) first on each vector $(f_{(j,k,\ell)})_{j=0}^{m-1}$, so that
$(f^{(j)}_{(k,\ell)})_{j=0}^{m-1} = \mathrm{DCT}(f_{(j,k,\ell)})_{j=0}^{m-1},\; \mathrm{for\;all}\;k,\ell.$
Then, perform another DCT, but this time on the vector $(f^{(j)}_{(k,\ell)})_{k=0}^{n-1}$, so that
$(f^{(j,k)}_{(\ell)})_{k=0}^{n-1} = \mathrm{DCT}(f^{(j)}_{(k,\ell)})_{k=0}^{n-1},\; \mathrm{for\;all}\;j,\ell.$
Finally, perform a discrete Fourier transform (DFT) on the $p$-tuple $(f^{(j,k)}_{(\ell)})_{\ell=0}^{p-1}$, so that
$(f^{(j,k,\ell)})_{\ell=0}^{p-1} = \mathrm{DFT}(f^{(j,k)}_{(\ell)})_{\ell=0}^{p-1},\; \mathrm{for\;all}\;j,k.$
Each value $f^{(j,k,\ell)}$ is then the coefficient on the basis function $q_{(j,k,\ell)}$. This series of transforms can be performed (or inverted) in $\mathcal{O}(N\log N)$ operations with the fast Fourier transform, where $N=mnp$.

Notably, it has been shown that this DCT-DCT-DFT transform offers a significant speed-up over comparable Legendre and Chebyshev-Legendre transforms, used with different polynomial bases \cite{shen_2001}. While this is of little concern for the heat equation itself, where it need only be done at the beginning and end of simulation, this becomes more significant when considering the nonlinear forcing we discuss in Section \ref{sec:navierstokes}. In the latter case, we need to perform a DCT-DCT-DFT transform once at each timestep, and the benefits of this fast transform are more significant.

\section{The Heat Equation}
The heat equation in $\mathit{C}$ is as follows:
$$
\left(\partial_t- \alpha\nabla^2\right)T = g(r,z,\theta,t),
$$
where $T(r,z,\theta,t)$ denotes the temperature at a given point in space and time in the cylinder, $g(r,z,\theta,t)$ is the rate of heat generation at a point $(r,z,\theta)$ at time $t$, and $\alpha>0$ is the ``thermal diffusivity'' of the domain, or the rate at which heat transfers to colder areas.

To integrate the heat equation forward in time, we use a high order \emph{backwards differentiation formula} (BDF), introduced in \cite{10.2307/88864}. This allows for a high degree of accuracy in time, while maintaining stability; in particular, as the heat equation is \emph{stiff}, implicit methods (such as BDF) are necessary for its efficient solution. The BDFs are given below for $b = 1,4$ with time-step $h$:
\bigskip

\paragraph{BDF 1}
$$
(1-h\partial_t)f(t+h) \approx f(t),
$$
\medskip

\paragraph{BDF 4}
$$
\left(1-\tfrac{12}{25}h\partial_t\right)f(t+h) \approx \tfrac{48}{25}f(t) - \tfrac{36}{25}f(t-h) + \tfrac{16}{25}f(t-2h) - \tfrac{3}{25}f(t-3h).
$$
\smallskip

By computing the right-hand side of one of the above expressions and inverting the differential operator on the left, we obtain an approximation of $f$ at the next timestep: to first order with BDF 1 and to fourth order with BDF 4.

More precisely, define $\delta^{(b)} f(t)$ to be the RHS of the $b^{th}$ order BDF, and define $\kappa^{(b)}$ to be the (positive) coefficient of $h\partial_tf(t+h)$ in the same equation. For instance, $\kappa^{(4)}=\frac{12}{25}h$ and $\delta^{(4)}f(t) =\frac{48}{25}f(t) - \frac{36}{25}f(t-h) + \frac{16}{25}f(t-2h) - \frac{3}{25}f(t-3h)$.
We can work this approximation into the heat equation by solving the latter for $\frac{\partial}{\partial t}T$:
\begin{equation}\label{heatBDF}
	\partial_tT(r,z,\theta,t)\approx \alpha\nabla^2 T(r,z,\theta,t) + g(r,z,\theta,t).
\end{equation}
Once we know $\delta^{(b)} f(t)$, which we would if we were progressing through equidistant time steps, we can plug (\ref{heatBDF}) into the $b^{th}$ order BDF and get a \emph{Helmholtz equation}:
$$
\left(1- \kappa^{(b)}\alpha\nabla^2\right) T(r,z,\theta,t+h) \approx \delta^{(b)} T(r,z,\theta,t) + \kappa^{(b)}g(r,z,\theta,t+h),
$$
where the RHS is already known and $I$ is the identity operator. An initial problem is in the apparent unboundedness of the operator on the LHS at $r=0$, as it includes the terms $r^{-1}$ and $r^{-2}$:
$$
1- \kappa^{(b)}\alpha\nabla^2 = 1 - \kappa^{(b)}\alpha\left( r^{-1}\partial_r+\partial_r^2 + r^{-2}\partial_\theta^2 + \partial_z^2\right).
$$
To resolve this, we multiply both sides of the Helmholtz equation by $r^2$ to obtain
\begin{equation}\label{Helmholtz}
r^2\left(I- \kappa^{(b)}\alpha\nabla^2\right) T(r,z,\theta,t+h) \approx r^2\delta^{(b)} T(r,z,\theta,t) + \kappa^{(b)}r^2g(r,z,\theta,t+h).
\end{equation}
The CCF$_{(m,n,p)}$ polynomial basis provides significant benefits for solving the above equation. If we work in ``coefficient space'', or the set of coefficients $(T^{(j,k,\ell)})$, the differential operator on the left becomes a sparse linear operator; as we will see, we can efficiently solve the resulting equation using an iterative procedure.

That is, we can ``discretize'' the left-hand operator in (\ref{Helmholtz}) in the space $P_{(m,n,p)}$. To start with, the $m\times n$ matrix $\mathbf{X}_\ell(t)$ is defined by $[\mathbf{X}_\ell(t)]_{j,k}=T^{(j,k,\ell)}(t)$, where $[\mathbf{X}_\ell(t)]_{j,k}$ is the element in the $j^{th}$ row and $k^{th}$ column of $\mathbf{X}_\ell(t)$. We similarly define the $m\times n$ matrix $\mathbf{g}_\ell(t)$ as the matrix of interpolation coefficients of $g(r,z,\theta,t)$.

From here, we discretize necessary operators as matrices to act on $\mathbf{X}_\ell(t)$.
\begin{itemize}
	\item The matrix $C_{01}$ converts a $T_i$ coefficient vector to a $U_i$ coefficient vector, and the matrix $C_{12}$ converts a $U_i \equiv C^{(1)}_i$ coefficient vector to a $C^{(2)}_i$ coefficient vector:
	$$C_{01} = 
	\begin{bmatrix*}[r] 0 & \frac{1}{2} & & &  \phantom{\ddots}\\
	\frac{1}{2} & 0 & \frac{1}{2} & &  \phantom{\ddots}\\
	& \frac{1}{2} & \ddots & \ddots& \phantom{\ddots}\\
	& & \ddots & \ddots &\frac{1}{2} \\
	\phantom{\ddots} & \phantom{\ddots} & & \frac{1}{2} & 0\end{bmatrix*},
	\qquad
	C_{12} = 
	\begin{bmatrix*}[r] 
	1 & 0 & -\frac{1}{3} & & \phantom{\ddots} \\
	& \frac{1}{2} & 0 & \ddots&  \phantom{\ddots} \\
	& & \frac{1}{3} & \ddots & -\frac{1}{n}\\
	& & & \ddots & 0\\
	\phantom{\ddots} & \phantom{\ddots} & \phantom{\ddots} & \phantom{\ddots} & \frac{1}{n}
	\end{bmatrix*}.
	$$
	The matrix $C_{02}$ is defined by $C_{02}=C_{12}C_{01}$. It converts a $T_i$ coefficient vector to a $C^{(2)}_i$ coefficient vector.
	
	\item The $b^{th}$-order BDF operator $\delta^{(b)}$ is found by sampling its argument at $(t, t-h, ..., t-(b-1)h)$ and taking differences, as with the RHS of the $b^{th}$ order BDF. For instance,
	$$\delta^{(4)}\mathbf{X}_\ell(t)=\frac{48}{25}\mathbf{X}_\ell(t) - \frac{36}{25}\mathbf{X}_\ell(t-h) + \frac{16}{25}\mathbf{X}_\ell(t-2h) - \frac{3}{25}\mathbf{X}_\ell(t-3h).$$
	
	\item The multiplication operator $R$ multiplies a $C^{(2)}_k(r)$ coefficient vector by $r$:
	$$
	R = 
	\begin{bmatrix*}[r] 
	0 & \frac{2}{3} & & & \phantom{\ddots}\\
	\frac{1}{4} & 0 & \frac{5}{8} & & \phantom{\ddots}\\
	& \frac{1}{3} & 0 & \ddots & \phantom{\ddots}\\
	& & \ddots & \ddots & \frac{n+2}{2n+2} \\
	\phantom{\ddots} & \phantom{\ddots} & \phantom{\ddots} & \frac{n-1}{2n} & 0 \\
	\end{bmatrix*}.
	$$
	Note that the subdiagonal pattern starts with $\frac{2-1}{4}=\frac{1}{4}$, not with $\frac{1-1}{2}=0$. For convenience, also define $R_2=R^2C_{02}$.

	\item The differentiation operator $D_1$ is defined as $D_1=C_{01}\partial_r$, where $\partial_r$ differentiates a $T_k(r)$ coefficient vector by $r$; similarly, the second differentiation operator is $D_2=C_{02}\partial_r^2$, where $\partial_r^2$ differentiates a $T_k(r)$ coefficient vector twice by $r$: 
	$$
	D_1 = 
	\begin{bmatrix*}[r] 
	0 & 1 & & & \phantom{\ddots}\\
	& & 2 & & \phantom{\ddots}\\
	& & & \ddots & \phantom{\ddots}\\
	& & & & n-1 \\
	\phantom{\ddots} & \phantom{\ddots} & \phantom{\ddots} & \phantom{\ddots} &0
	\end{bmatrix*},
	\qquad
	D_2 = 
	\begin{bmatrix*}[r] 
	0 & 0 & 4 & & \phantom{\ddots} \\
	& & & \ddots & \phantom{\ddots} \\
	\phantom{\ddots} & & & & 2n-2 \\
	\phantom{\ddots} & & & & 0\\
	\phantom{\ddots} &\phantom{\ddots} &\phantom{\ddots} &\phantom{\ddots} & 0\\
	\end{bmatrix*}.
	$$
\end{itemize}
With these operators, (\ref{Helmholtz}) can be discretized (with the extrapolation $g(t+h)\approx g(t)$) as the following matrix equation:
\begin{align}\label{HelmholtzMatEqn}
	\left(R_2-\kappa^{(b)}\alpha(RC_{12}D_1+R^2D_2-\ell^2C_{02})\right)\mathbf{X}_\ell(t+h)C^T_{02}-\kappa^{(b)}\alpha R_2\mathbf{X}_\ell(t+h)D^T_2 = \nonumber\\
	R_2\left(\delta^{(b)} \mathbf{X}_\ell(t)+\kappa^{(b)}\mathbf{g}_\ell(t)\right)&C^T_{02}.
\end{align}
The discretization process used above is discussed further in \cite{doi:10.1137/120865458}. There are two notable properties of (\ref{HelmholtzMatEqn}): each value of $\ell$ can be solved for independently (we can \textquotedblleft decouple in $\theta$\textquotedblright{}), and every matrix used in the construction of the equation is very sparse. We use the ADI method to solve equation (\ref{HelmholtzMatEqn}) in $\mathcal{O}(N\log N)$ operations \cite{doi:10.1137/0103003}.

\begin{algo}[ADI Method]\label{alg1}
	Let $A,C$ be $m\times m$ complex matrices, $B,D$ be $n\times n$ complex matrices, and $X, E$ be $m\times n$ complex matrices. Set two complex-valued sequences $(u_j)_{j=0}^{J-1}$ and $(v_j)_{j=0}^{J-1}$, and set $X_0 =0$.
	
	The \emph{ADI method} for approximately solving the equation $AXB+CXD=E$ for $X$ is as follows:
	\begin{enumerate}
		\item Solve $(u_jC-A)X_{j+\frac{1}{2}}B=CX_j(u_jB+D)+E$ for $X_{j+\frac{1}{2}}$.
		\item Solve $CX_{k+1}(v_jB+D)= (v_jC-A)X_{j+\frac{1}{2}}B-E$ for $X_{j+1}$.
		\item Repeat steps 1 and 2 for $J$ iterations; $X_J \approx X$.
	\end{enumerate}
	With an appropriate choice of $(u_j)_{j=0}^{J-1}$ and $(v_j)_{j=0}^{J-1}$, we can get sufficient accuracy with $J=\mathcal{O}(\log mn)$. This iteration then gives an accurate solution in $\mathcal{O}(mn\log mn)$ operations, assuming that the initial Sylvester equation is well-conditioned. All such equations in this work are well-conditioned unless otherwise stated.
\end{algo}

To ensure quasi-optimal convergence, it is critical to choose shifts $u_k$ and $v_k$ appropriately. This can be done using \cite[Theorem 2.1]{10.1093/imanum/drz034}, following their logic to ensure the requisite hypotheses on the matrices $A$, $B$, $C$, and $D$. For this, first define $a$ and $b$ to be the minimum and maximum eigenvalues of $C^{-1}A$ and similarly for $c$ and $d$ and $-DB^{-1}$. Set $\gamma = |c-a||d-b|/(|c-b||d-a|)$ and $\alpha = -1+2\gamma +2\sqrt{\gamma^2-\gamma}$, and finally $k = \sqrt{1-1/\alpha^2}$. Then in the setting of Algorithm \ref{alg1}, we find optimal values of
\[u_j = T\left(-\alpha\text{dn}\left[\frac{2j+1}{2J}K(k),k\right]\right),\qquad v_j=T\left(\alpha\text{dn}\left[\frac{2j+1}{2J}K(k),k\right]\right),\]
where $K(k)$ is the complete elliptic integral of the first kind and $\operatorname{dn}(z,k)$ is the Jacobi elliptic function of the third kind (see \cite{822801} for details), and $T$ is the Möbius transformation mapping $\{-\alpha,-1,1,\alpha\}$ to $\{a,b,c,d\}$. It is shown in \cite{10.1093/imanum/drz034} that this choice of shift parameters leads to an exponential convergence in the $L^2$ norm.

It should be noted that this does not apply to arbitrary $A$, $B$, $C$, and $D$. Rather, following the logic of \cite{10.1093/imanum/drz034}, we can transform our equation into a standard Sylvester equation $QX - XP = F$, where $Q$ and $P$ are normal and with disjoint, real spectra. In turn, this demonstrates the convergence with the above parameters.

\begin{figure}
	\begin{subfigure}[b]{.33\textwidth}
		\centering
		\includegraphics[scale=.6]{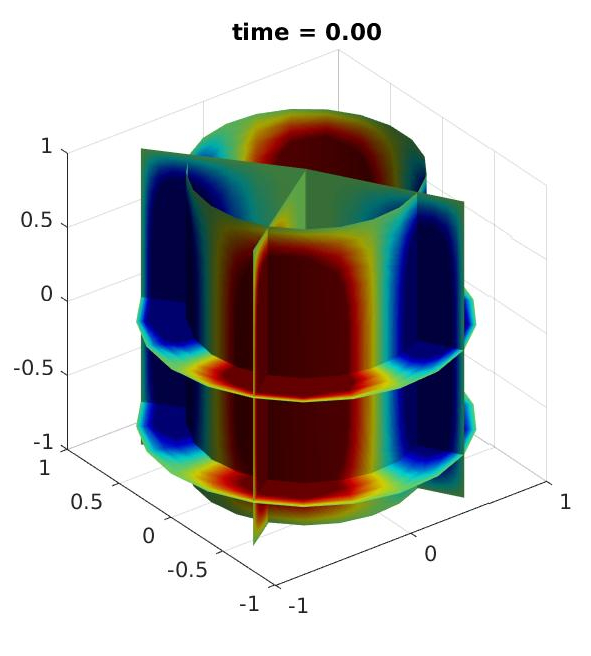}
	\end{subfigure}%
	\begin{subfigure}[b]{.33\textwidth}
		\centering
		\includegraphics[scale=.6]{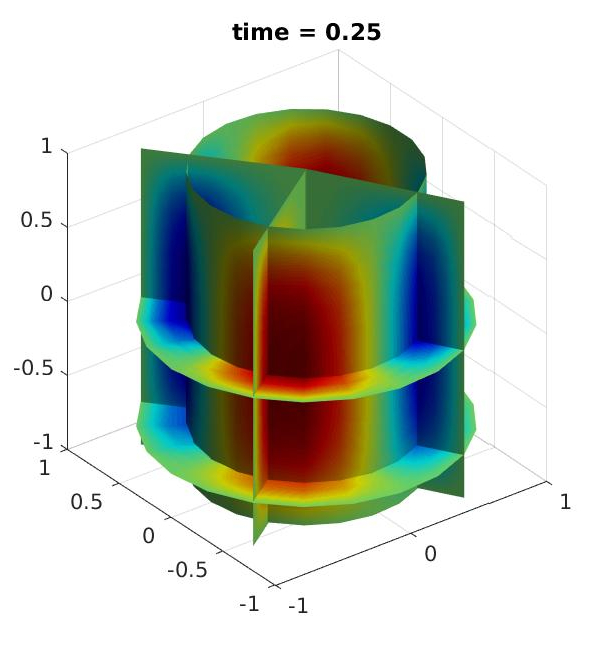}
	\end{subfigure}%
	\begin{subfigure}[b]{.33\textwidth}
		\centering
		\includegraphics[scale=.6]{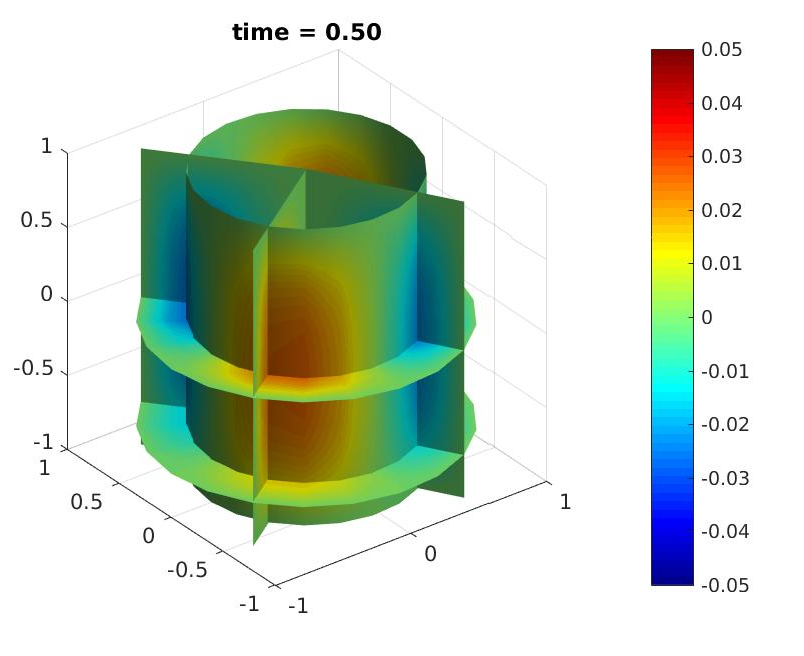}
	\end{subfigure}
	\caption{A simulation of the heat equation in the cylinder $\mathit{C}$. This is shown on various slices of the cylinder, where the temperature value is given by the color bars on the right of each figure. There is no advective motion, making the heat equation computationally faster to solve than the Navier--Stokes equations.}\label{heatfig}
\end{figure}

In short, at each time step, we can solve for each $\mathbf{X}_\ell$ in $\mathcal{O}(mn\log mn)$ operations. There are $p$ choices of $\ell$, though, so this is scaled up to $\mathcal{O}(mnp\log mn)$ total operations at each time step. Since $N = mnp$, the final algorithm has a total of  $\mathcal{O}(N\log N)$ operations. This compares to times of $\mathcal{O}(N^{7/3})$, or $\mathcal{O}(m^3n^3p)$, with current spectral methods. A solution example is illustrated in Figure \ref{heatfig}.

\section{Application to the Incompressible Navier--Stokes Equations}\label{sec:navierstokes}
In the following section, we highlight an application of the preceding method to solve the incompressible Navier--Stokes equations in the cylinder. In vorticity form (i.e., taking the curl of its primitive form), the equations are written as
\begin{equation}\label{eq:ns}
\left(\partial_t -\tfrac{1}{\text{Re}}\nabla^2\right)\vec{\omega} = \nabla\times[\vec{v}\times\vec{\omega}],
\end{equation}
where $\vec{\omega}:=\nabla\times\vec{v}$ is the \emph{vorticity} of the flow.

\begin{figure}
	\centering
	\begin{subfigure}{.26\textwidth}
		\centering
		\raisebox{-.5\height}{ \includegraphics[scale=.24]{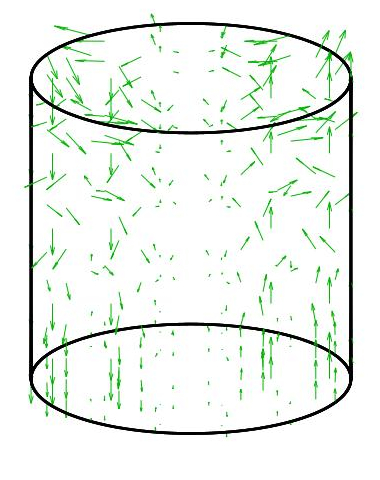}}
		\caption*{Vorticity Field}
	\end{subfigure}
	{\huge $=$}
	\begin{subfigure}{.26\textwidth}
		\centering
		\raisebox{-.5\height}{ \includegraphics[scale=.24]{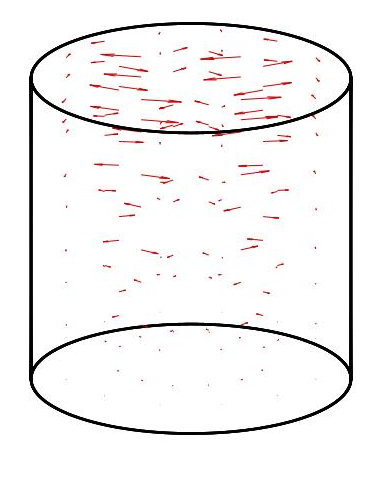}}
		\caption*{Toroidal Field}
	\end{subfigure}
	{\huge $+$}
	\begin{subfigure}{.26\textwidth}
		\centering
		\raisebox{-.5\height}{ \includegraphics[scale=.24]{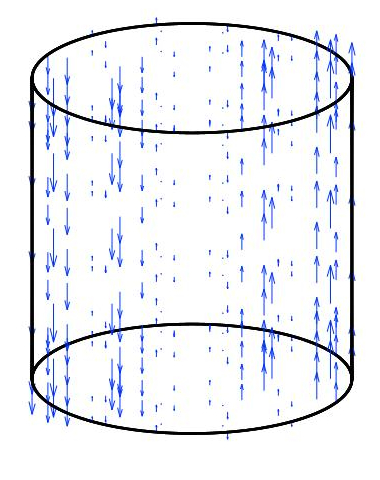}}
		\caption*{Poloidal Field}
	\end{subfigure}
	
	\caption{An example of a vorticity vector field being decomposed into its poloidal and toroidal components. The components are each given by a single scalar field. The poloidal and toroidal fields were computed numerically.}\label{PTfig}
\end{figure}

\paragraph{The Poloidal--Toroidal Decomposition} To solve the Navier--Stokes equations, we introduce the poloidal--toroidal (PT) decomposition to split incompressible vector fields into two scalar components each. Specifically, for any incompressible vector field $\vec{V}:\mathit{C}\rightarrow \mathbb{R}^3$ on the cylinder, we can find two scalar fields $\lambda$ and $\gamma$ such that
\begin{equation}\label{eq:pt}
	\vec{V} = \nabla\times[\lambda\hat{z}]+\nabla\times\nabla\times[\gamma\hat{z}].
\end{equation}
These fields are called the \emph{toroidal} and \emph{poloidal} components of $\vec{V}$, respectively. The new variables allow for parallelization and maintain the incompressibility of the fluid automatically. Incompressibility is traditionally handled using projection methods---i.e., by projecting the solution onto the space of incompressible vector fields after each time step---but these methods can introduce large errors in the simulation \cite{10.2307/2153791}. An example PT decomposition is shown in Fig. \ref{PTfig}.

We introduce here only a few critical properties of the PT decomposition, used in our proposed algorithm. For a more general discussion of the PT decomposition, please see \cite{Boronski_2007}.

For one, the PT decomposition of a field $\vec{V}$ can be found by solving a disk Poisson equation. In the above notation, a short calculation shows that
\begin{equation}\label{eq:ptpoisson}
	\hat{z} \cdot \vec{V} = -\Delta_h \gamma,
	\qquad\quad
	\hat{z} \cdot \nabla\times \vec{V} = -\Delta_h\lambda,
\end{equation}
where $\Delta_h = \nabla^2-\partial_z^2$ is the horizontal Laplacian. 

Secondly, if $\vec{W}=\nabla\times\vec{V}$, then we can write
\begin{equation}\label{eq:ptcurl}
	\vec{W} = \nabla\times[-\nabla^2\gamma\hat{z}]+\nabla\times\nabla\times[\lambda\hat{z}].
\end{equation}
Conversely, this means that, given a PT decomposition $(\lambda_W,\gamma_W)$ of $\vec{W}$, we can calculate $(\lambda,\gamma)$ using
\[\nabla^2\gamma = -\lambda_W,\qquad\quad \lambda = \gamma_W.\]
This Poisson solve can be performed using the methods developed above.

\paragraph{Rewriting the Navier--Stokes Equations} Suppose we have found PT decompositions $(\lambda_\omega,\gamma_\omega)$ of $\vec{\omega}:=\nabla\times\vec{v}$ and $(\lambda_a,\gamma_a)$ of $\nabla\times(\vec{v}\times\vec{\omega})$, for instance, by solving (\ref{eq:ptpoisson}). Introducing the decomposition (\ref{eq:pt}) of these two fields into (\ref{eq:ns}) and splitting by poloidal and toroidal components, we find the two equations
\begin{equation}\label{eq:NSrewrite}
	\left(\partial_t - \tfrac{1}{Re} \nabla^2 \right)\lambda_\omega = \lambda_a, \qquad \quad \left(\partial_t - \tfrac{1}{Re} \nabla^2 \right)\gamma_\omega = \gamma_a.
\end{equation}
The proposed method then proceeds as follows. Starting with the PT decomposition of the vorticity $\vec{\omega}=\nabla\times\vec{v}$, calculate the decomposition $(\lambda_v,\gamma_v)$ of the velocity $\vec{v}$ by solving (\ref{eq:ptcurl}). Then, calculate $\nabla\times[\vec{v}\times\vec{\omega}]$ explicitly and use (\ref{eq:ptpoisson}) to find $(\lambda_a,\gamma_a)$ at the current time step. 

At this point, the backwards differentiation formulas we use for the heat equation are no longer feasible for a nonlinear equation; each time step would require an expensive rootfinding algorithm in order to invert the now-nonlinear operator on the lefthand side of the BDF. Instead, we apply implicit--explicit versions of these formulas, following \cite{10.2307/2158449}. These integrate the diffusive, linear component as before, but now integrate the nonlinear components of the Navier--Stokes equations explicitly.

All in all, this requires three applications of our ADI-based Poisson solver at each timestep: one to find $\gamma_v$ at the current timestep, and two to find $(\lambda_\omega,\gamma_\omega)$ at the next timestep. Further, it requires solving two ``horizontal'' Poisson equations, given by (\ref{eq:ptpoisson}), which can be done in optimal time by diagonalizing in both $\theta$ and $z$ directions. Finally, the explicit computation of $\nabla\times[\vec{v}\times\vec{\omega}]$ can be done in quasi-optimal time; though it requires transforming $\vec{v}$ and $\vec{\omega}$ back to value space to diagonalize the multiplication operator, these transforms can be done in $\mathcal{O}(N\log N)$ time themselves.
\section{Numerical Experiments}\label{sec:num}
Three methods for solving the forced heat equation were tested against an exact solution using various discretization sizes. Namely, we compared our new algorithm against a spectral collocation method on the same CCF grid (see Definition \ref{grid}) and against a finite difference method on an equidistant grid in $(r,z,\theta)$ space. We summarize our results in Table \ref{table}, including the maximum relative error values for these methods as well as the computational times required to find these solutions. 

Our new spectral method and an implemented spectral collocation method were run for 200 time steps on the same grid. The error is comparable in the case $N=7^3$ and $N=15^3$, though the latter shows substantially higher error in the $N=11^3$ and $N=19^3$ cases, possibly due to ill-conditioning of the applied matrices. Nonetheless, the asymptotic difference in runtime becomes clear for $N = 11^3, 15^3, 19^3$. 

We implemented a basic finite difference method for comparison, using second-order-accurate differentiation operators and an LU decomposition of the resulting sparse system. Note that finite difference methods can be optimized further than this; for one, our application of the ADI algorithm could be developed in a finite-difference case as well. This method was not able to achieve the same accuracy on any feasible mesh size, so we instead record their error and runtime in finer meshes of size $N = 61^3$, $81^3$, $101^3$, and $121^3$. It was also run for only 6 time steps in each case, as compared to the 200 tested with spectral methods, so we expect that the total runtime would be $\sim 200/6$ of the reported values.

Overall, these results suggest that our new spectral method represents a steep improvement in accuracy over non-spectral methods and a steep improvement in runtime over competing spectral methods.


\begin{table}
	\footnotesize
	\centering
	\begin{tabular}{ l cccc}
		\toprule\multicolumn{1}{c}{} & \multicolumn{4}{c}{Maximum Error} \\\cmidrule(lr){2-5}
		Method \textbackslash{} $N$ & $7^3$ & $11^3$ & $15^3$ & $19^3$ \\ \midrule
		New Spectral Method & $6.03\cdot 10^{-5}$ & $5.88\cdot10^{-5}$ & $3.22\cdot10^{-5}$ & $1.76\cdot10^{-5}$ \\
		Spectral Collocation & $6.40\cdot10^{-5}$ & $5.85\cdot10^{-3}$ & $3.51\cdot10^{-5}$ & $5.90\cdot10^{-3}$ \\
		Finite Difference & $>3.80\cdot10^{-1}$ & $>2.93\cdot10^{-1}$ & $>2.39\cdot10^{-1}$ & $>2.00\cdot10^{-1}$ \\\bottomrule
		\multicolumn{5}{c}{}\\\multicolumn{5}{c}{}\\\toprule
		\multicolumn{1}{c}{} & \multicolumn{4}{c}{Computational Time (sec)}\\\cmidrule(lr){2-5}
		Method \textbackslash{} $N$ & $7^3$ & $11^3$ & $15^3$ & $19^3$ \\ \midrule
		New Spectral Method & $1.00$ & $3.19$ & $5.88$ & $6.50$\\
		Spectral Collocation& $1.06$ & $4.42$ & $11.9$ & $17.9$ \\
		Finite Difference  & $>25.1$ & $>28.3$ & $>56.8$ & $>100$ \\\bottomrule
	\end{tabular}
	\caption{Direct comparison of our algorithm against a spectral collocation method and a finite difference method for the heat equation. Above are the maximum relative error values for these methods (measured against an exact solution), as well as the computational times required to find these solutions. The first two methods were run for 200 time steps, with $N = 7^3$, $11^3$, $15^3$, and $19^3$ respectively. Due to machine limitations, the finite difference method was only run for 6 time steps, with $N = 61^3$, $81^3$, $101^3$, and $121^3$ respectively; values were marked with ``$>$'' to denote this difference. The new spectral method achieves the highest accuracy in the least computational time of the tested methods.}\label{table}
\end{table}

\section{Conclusion}
There are several published spectral methods which can be used to approximate solutions to the heat equation in a closed cylinder \cite{gottlieb1977numerical, TAN198581, PULICANI199193, PRIYMAK1995366, doi:10.1137/S1064827595295301,CiCP-5-426}, but to our knowledge, with the exception of \cite{doi:10.1137/S1064827595295301}, they do not achieve quasi-optimal runtimes. The latter is able to achieve this computational complexity using Legendre or Chebyshev-Legendre polynomial bases, but (a) their approach requires a boundary condition at the central axis of the cylinder and (b) it has been found that the corresponding coefficient transforms of these polynomials are far more expensive than the DCT-DCT-DFT transform discussed in Section \ref{sec:approx} \cite{shen_2001}. Note that the latter is of minor concern when solving the heat equation itself, but becomes critical when applying these methods to systems with nonlinear forcing---as discussed in Section \ref{sec:navierstokes}, a coefficient transform must be performed at each time step in solving the Navier--Stokes equations. Using Chebyshev polynomials, their method achieves a runtime of $\mathcal{O}(N^{4/3})$.

The current method resolves this problem by leveraging the ADI algorithm of \cite{doi:10.1137/0103003}. In short, discretizing along a Chebyshev-Chebyshev-Fourier (CCF) grid in the cylinder allows us to decompose the heat or Poisson equation into a collection of sparse Sylvester equations, while resolving all parts of the cylinder with high accuracy and allowing the use of a fast coefficient transform to recover function values. These Sylvester equations can be solved by the iterative procedure of Algorithm \ref{alg1}; by choosing appropriate algorithm parameters, we can guarantee convergence in $\mathcal{O}(\log N)$ iterations, for a total of $\mathcal{O}(N\log N)$ operations. As demonstrated in Section \ref{sec:num}, this accuracy generally represents an improvement of several orders of magnitude over finite difference methods, and the quasi-optimal runtime is several times faster than spectral collocation methods in even small system sizes ($N\leq 19^3$).
\section*{Appendix: Boundary Conditions for the Navier--Stokes Equations}
The formulation (\ref{eq:NSrewrite}) of the Navier--Stokes equations requires three sets of boundary conditions. First, conditions for computing the PT decomposition of a given vector field. Second, conditions for finding the poloidal component of $\vec{v}$ from the toroidal component of $\vec{\omega}$, as in (\ref{eq:ptcurl}). Finally, conditions for solving the discretized Navier--Stokes equations (\ref{eq:NSrewrite}) for $\vec{\omega}$.

For the PT decomposition $\vec{V} = \nabla\times[\lambda\hat{z}]+\nabla\times\nabla\times[\gamma\hat{z}]$ to be correct, it must give correct results for each of the three vector components $(V^r\hat{r},V^\theta\hat{\theta},V^z\hat{z})$ of $\vec{V}$. Since the $\hat{z}$ component is folded into the horizontal Poisson equation, we are left with $\partial_\theta\lambda + r\partial_r\partial_z\gamma = rV^r$ and $\partial_\theta\partial_z\gamma - r\partial_r\lambda = rV^\theta$. These equations only need to be enforced on the boundary of the cylinder, and can be discretized in the same manner as were the Helmholtz and Poisson equations.

In order to find boundary conditions for the poloidal component of $\vec{v}$, we think instead in terms of the vector potential $\vec{\Psi}$, defined by $\nabla\times\vec{\Psi}=\vec{v}$. Suppose $\vec{\Psi}$ has a decomposition $(\lambda_\psi,\gamma_\psi)$. From (\ref{eq:ptcurl}), we see that $\lambda_\psi = \gamma_v$ satisfies $\nabla^2\lambda_\psi= -\lambda_\omega$ and that $\nabla^2\gamma_\psi=-\gamma_\omega$, and we propose solving both simultaneously. This leads to one more Poisson solve per timestep, but does not change the quasi-optimality of the solution. For \emph{no-slip} boundary conditions on $\lambda_\psi = \gamma_v$---that is, where the boundary of the cylinder is non-moving and frictional---the vector potential must be normal to the boundary of the cylinder (which can be seen from Stokes' Theorem). When $z = \pm 1$, the $\hat{r}$ and $\hat{\theta}$ components of $\Psi$ must vanish: $\partial_\theta\lambda_\psi + r\partial_r\partial_z\gamma_\psi = \partial_\theta\partial_z\gamma_\psi - r\partial_r\lambda_\psi = 0$. When $r = 1$, the $\hat{z}$ and $\hat{\theta}$ components of $\Psi$ must vanish: $(r\partial_r + r^2\partial^2_r +\partial^2_\theta)\gamma_\psi = \partial_\theta\partial_z\gamma_\psi - r\partial_r\lambda_\psi = 0$. These equations can be discretized as before.

For solving the Navier--Stokes equations, we continue to use no-slip conditions. Due to work by Quartepelle \cite{https://doi.org/10.1002/fld.1650010204}, it is known that this no-slip condition on the velocity field is equivalent to the vorticity field being $L^2$ orthogonal to all incompressible vector fields $\vec{\eta}$ with $\nabla^2\vec{\eta} = 0$. To form these $\vec{\eta}$, first let $(\eta_i)$ be a basis of scalar cylindrical harmonics. For some $i$, let $\vec{\eta}_i=\nabla\times[\eta_i \hat{z}]$ and $\vec{\eta}_i^{(2)} = \nabla\times\nabla\times[\eta_j \hat{z}].$ As all incompressible vector fields have a PT decomposition and the Laplacian commutes with the curl, the set $\{\vec{\eta}_i\}\cup\{\vec{\eta}_i^{(2)}\}$ spans the set of $\vec{\eta}$. To achieve $L^2$ orthogonality between $\vec{\omega}$ and $\vec{\eta}$ in $C$, the following equation must be met:
\begin{align*}
\int_C\big(\vec{\eta}\cdot(\nabla\times[\lambda_\omega \hat{z}]+&\nabla\times\nabla\times[\gamma_\omega \hat{z}])\big)\;dV \\&= 
\int_C\left([\nabla\times\vec{\eta}]\cdot \lambda_\omega \hat{z}\right)dV+\int_{\partial C}\left(\vec{n}\cdot[\lambda_\omega \hat{z}\times\vec{\eta}]\right)dS\\&
\quad+\int_{\partial C}\left(\vec{n}\cdot[\gamma_\omega \hat{z}\times\nabla\times\vec{\eta}
+(\nabla\times\vec{\eta})\times\gamma_\omega \hat{z}]\right)dS=0.
\end{align*}
Let $\vec{\eta}=\eta^r \hat{r} + \eta^\theta \hat{\theta} + \eta^z \hat{z}$. Then the above equation simplifies to
\begin{align*}
\int\limits_{r=1}\eta^\theta\lambda_\omega\;d\theta dz &+ \int_C[\eta^\theta+r\partial_r\eta^\theta-\partial_\theta\eta^r]\lambda_\omega\;drd\theta dz \\&+ 
\int\limits_{r=1}\big([\partial_z\eta^r-\partial_r\eta^z]\gamma_\omega +
\eta^z\partial_r\gamma_\omega\big)\;d\theta dz \\&- \int\limits_{z=1}\left(\eta^\theta\partial_\theta\gamma_\omega + r\eta^r\partial_r\gamma_\omega\right)\;drd\theta \\&+
\int\limits_{z=-1}\left(\eta^\theta\partial_\theta\gamma_\omega + r\eta^r\partial_r\gamma_\omega\right)\;drd\theta= 0.
\end{align*}
This equation can be discretized as before, creating one boundary condition for each $\vec{\eta}$. These conditions can also be included without \textquotedblleft coupling in $\theta$\textquotedblright{} by first enforcing a simpler set of boundary conditions and then performing a low-rank update on the solution matrix.
\pagebreak

\bibliographystyle{siam}
\bibliography{bibthing}

\begin{thebibliography}{10}

\bibitem{abramowitz1965handbook}
{\sc M.~Abramowitz and I.~Stegun}, {\em Handbook of Mathematical Functions:
  With Formulas, Graphs, and Mathematical Tables}, Applied mathematics series,
  Dover Publications, 1965.

\bibitem{10.2307/2158449}
{\sc U.~M. Ascher, S.~J. Ruuth, and B.~T.~R. Wetton}, {\em Implicit-explicit
  methods for time-dependent partial differential equations}, SIAM Journal on
  Numerical Analysis, 32 (1995), pp.~797--823.

\bibitem{CiCP-5-426}
{\sc F.~Auteri and L.~Quartapelle}, {\em Spectral elliptic solvers in a finite
  cylinder}, Communications in Computational Physics, 5 (2009), pp.~426--441.

\bibitem{doi:10.1137/080737058}
{\sc S.~Berrone}, {\em A local-in-space-timestep approach to a finite element
  discretization of the heat equation with a posteriori estimates}, SIAM
  Journal on Numerical Analysis, 47 (2009), pp.~3109--3138.

\bibitem{Boronski_2007}
{\sc P.~Boronski and L.~S. Tuckerman}, {\em Poloidal{\textendash}toroidal
  decomposition in a finite cylinder. i: Influence matrices for the
  magnetohydrodynamic equations}, Journal of Computational Physics, 227 (2007),
  pp.~1523--1543.

\bibitem{Boull__2020}
{\sc N.~Boull{\'{e} } and A.~Townsend}, {\em Computing with functions in the
  ball}, {SIAM} Journal on Scientific Computing, 42 (2020), pp.~C169--C191.

\bibitem{https://doi.org/10.48550/arxiv.2103.16638}
{\sc N.~Boullé, J.~Słomka, and A.~Townsend}, {\em An optimal complexity
  spectral method for {Navier--Stokes} simulations in the ball}, 2021.

\bibitem{BOYD20111408}
{\sc J.~P. Boyd and F.~Yu}, {\em Comparing seven spectral methods for
  interpolation and for solving the poisson equation in a disk: Zernike
  polynomials, logan–shepp ridge polynomials, chebyshev–fourier series,
  cylindrical robert functions, bessel–fourier expansions, square-to-disk
  conformal mapping and radial basis functions}, Journal of Computational
  Physics, 230 (2011), pp.~1408--1438.

\bibitem{Burns_2020}
{\sc K.~J. Burns, G.~M. Vasil, J.~S. Oishi, D.~Lecoanet, and B.~P. Brown}, {\em
  Dedalus: A flexible framework for numerical simulations with spectral
  methods}, Physical Review Research, 2 (2020).

\bibitem{doi:10.1137/15M1035094}
{\sc M.~F. Causley, H.~Cho, A.~J. Christlieb, and D.~C. Seal}, {\em Method of
  lines transpose: High order l-stable {O(N)} schemes for parabolic equations
  using successive convolution}, SIAM Journal on Numerical Analysis, 54 (2016),
  pp.~1635--1652.

\bibitem{cooper2012introduction}
{\sc J.~M. Cooper}, {\em Introduction to Partial Differential Equations with
  MATLAB}, Applied and Numerical Harmonic Analysis, Birkh{\"a}user Boston,
  2012.

\bibitem{10.2307/88864}
{\sc C.~F. Curtiss and J.~O. Hirschfelder}, {\em Integration of stiff
  equations}, Proceedings of the National Academy of Sciences of the United
  States of America, 38 (1952), pp.~235--243.

\bibitem{mycode}
{\sc D.~Darrow}, {\em Spectral heat equation solver}.
\newblock \url{https://github.com/ddarrow90/spectral-heat-equation}, 2017.

\bibitem{doi:10.1137/0916061}
{\sc B.~Fornberg}, {\em A pseudospectral approach for polar and spherical
  geometries}, SIAM Journal on Scientific Computing, 16 (1995), pp.~1071--1081.

\bibitem{https://doi.org/10.48550/arxiv.1710.11259}
{\sc D.~Fortunato and A.~Townsend}, {\em Fast poisson solvers for spectral
  methods}, 2017.

\bibitem{10.1093/imanum/drz034}
{\sc D.~Fortunato and A.~Townsend}, {\em {Fast Poisson solvers for spectral
  methods}}, IMA Journal of Numerical Analysis, 40 (2019), pp.~1994--2018.

\bibitem{GODON1993171}
{\sc P.~Godon and G.~Shaviv}, {\em A two-dimensional time dependent chebyshev
  method of collocation for the study of astrophysical flows}, Computer Methods
  in Applied Mechanics and Engineering, 110 (1993), pp.~171--194.

\bibitem{gottlieb1977numerical}
{\sc D.~Gottlieb and S.~A. Orszag}, {\em Numerical analysis of spectral
  methods: theory and applications}, SIAM, 1977.

\bibitem{Krauthammer1979AccuracyOT}
{\sc T.~Krauthammer}, {\em Accuracy of the finite element method near a curved
  boundary}, Computers \& Structures, 10 (1979), pp.~921--929.

\bibitem{KWAN2009170}
{\sc Y.-Y. Kwan}, {\em Efficient spectral-galerkin methods for polar and
  cylindrical geometries}, Applied Numerical Mathematics, 59 (2009),
  pp.~170--186.

\bibitem{10.5555/1355322}
{\sc R.~LeVeque}, {\em Finite Difference Methods for Ordinary and Partial
  Differential Equations: Steady-State and Time-Dependent Problems (Classics in
  Applied Mathematics Classics in Applied Mathemat)}, Society for Industrial
  and Applied Mathematics, USA, 2007.

\bibitem{MATSUSHIMA1995365}
{\sc T.~Matsushima and P.~Marcus}, {\em A spectral method for polar
  coordinates}, Journal of Computational Physics, 120 (1995), pp.~365--374.

\bibitem{822801}
{\sc F.~Olver, D.~Lozier, R.~Boisvert, and C.~Clark}, {\em The NIST Handbook of
  Mathematical Functions}, Cambridge University Press, New York, NY, 2010-05-12
  00:05:00 2010.

\bibitem{doi:10.1137/120865458}
{\sc S.~Olver and A.~Townsend}, {\em A fast and well-conditioned spectral
  method}, SIAM Review, 55 (2013), pp.~462--489.

\bibitem{PhysRevLett.26.1100}
{\sc S.~A. Orszag}, {\em Galerkin approximations to flows within slabs,
  spheres, and cylinders}, Phys. Rev. Lett., 26 (1971), pp.~1100--1103.

\bibitem{doi:10.1137/0103003}
{\sc D.~W. Peaceman and H.~H. Rachford, Jr.}, {\em The numerical solution of
  parabolic and elliptic differential equations}, Journal of the Society for
  Industrial and Applied Mathematics, 3 (1955), pp.~28--41.

\bibitem{PRIYMAK1995366}
{\sc V.~Priymak}, {\em Pseudospectral algorithms for navier-stokes simulation
  of turbulent flows in cylindrical geometry with coordinate singularities},
  Journal of Computational Physics, 118 (1995), pp.~366--379.

\bibitem{PULICANI199193}
{\sc J.~Pulicani and J.~Ouazzani}, {\em A fourier-chebyshev pseudospectral
  method for solving steady 3-d navier-stokes and heat equations in cylindrical
  cavities}, Computers \& Fluids, 20 (1991), pp.~93--109.

\bibitem{https://doi.org/10.1002/fld.1650010204}
{\sc L.~Quartapelle and F.~Valz-Gris}, {\em Projection conditions on the
  vorticity in viscous incompressible flows}, International Journal for
  Numerical Methods in Fluids, 1 (1981), pp.~129--144.

\bibitem{10.2307/2153791}
{\sc J.~Shen}, {\em On error estimates of the projection methods for the
  {Navier--Stokes} equations: Second-order schemes}, Mathematics of
  Computation, 65 (1996), pp.~1039--1065.

\bibitem{doi:10.1137/S1064827595295301}
{\sc J.~Shen}, {\em Efficient spectral-galerkin methods iii: Polar and
  cylindrical geometries}, SIAM Journal on Scientific Computing, 18 (1997),
  pp.~1583--1604.

\bibitem{SHEN2000183}
{\sc J.~Shen}, {\em A new fast chebyshev–fourier algorithm for poisson-type
  equations in polar geometries}, Applied Numerical Mathematics, 33 (2000),
  pp.~183--190.

\bibitem{shen_2001}
{\sc J.~Shen}, {\em Efficient chebyshev-legendre galerkin methods for elliptic
  problems}, Houston J Math,  (2001).

\bibitem{https://doi.org/10.48550/arxiv.2111.04585}
{\sc C.~Strössner and D.~Kressner}, {\em Fast global spectral methods for
  three-dimensional partial differential equations}, 2021.

\bibitem{10.2307/2156231}
{\sc P.~N. Swarztrauber}, {\em A direct method for the discrete solution of
  separable elliptic equations}, SIAM Journal on Numerical Analysis, 11 (1974),
  pp.~1136--1150.

\bibitem{TAN198581}
{\sc C.~Tan}, {\em Accurate solution of three-dimensional poisson's equation in
  cylindrical coordinate by expansion in chebyshev polynomials}, Journal of
  Computational Physics, 59 (1985), pp.~81--95.

\bibitem{torres1999}
{\sc D.~Torres and E.~Coutsias}, {\em Pseudospectral solution of the
  two-dimensional {Navier--Stokes} equations in a disk}, SIAM Journal on
  Scientific Computing, 21 (1999), pp.~378--403.

\bibitem{10.5555/3384673}
{\sc L.~N. Trefethen}, {\em Approximation Theory and Approximation Practice,
  Extended Edition}, SIAM-Society for Industrial and Applied Mathematics,
  Philadelphia, PA, USA, 2019.

\bibitem{VASIL201653}
{\sc G.~M. Vasil, K.~J. Burns, D.~Lecoanet, S.~Olver, B.~P. Brown, and J.~S.
  Oishi}, {\em Tensor calculus in polar coordinates using jacobi polynomials},
  Journal of Computational Physics, 325 (2016), pp.~53--73.

\bibitem{verkley_97}
{\sc W.~Verkley}, {\em A spectral model for two-dimensional incompressible
  fluid flow in a circular basin i. mathematical formulation}, Journal of
  Computional Physics,  (1997).

\bibitem{doi:10.1137/16M1070207}
{\sc H.~Wilber, A.~Townsend, and G.~B. Wright}, {\em Computing with functions
  in spherical and polar geometries ii. the disk}, SIAM Journal on Scientific
  Computing, 39 (2017), pp.~C238--C262.

\bibitem{ZHIYIN201511}
{\sc Y.~Zhiyin}, {\em Large-eddy simulation: Past, present and the future},
  Chinese Journal of Aeronautics, 28 (2015), pp.~11--24.

\end{thebibliography}

\end{document}